\documentclass{amsart}
\usepackage{amsmath,amsthm,amssymb,array}
\usepackage[all]{xy}

\theoremstyle{plain}
\newtheorem{theorem}{Theorem}
\newtheorem{lemma}{Lemma}
\newtheorem{proposition}{Proposition}
\newtheorem{corollary}{Corollary}

\theoremstyle{definition}
\newtheorem{example}{Example}
\newtheorem{definition}{Definition}

\newcommand\Z{{\mathbb Z}}

\newcommand\Aut{{\mathsf{Aut}}}
\newcommand\Ker{{\mathsf{Ker}}}
\newcommand\Sym{{\mathsf{Sym}}}

\newcommand\G{{\mathcal G}}
\newcommand\A{{\mathcal A}}
\newcommand\F{{\mathcal F}}
\newcommand\dimh{{\dim_{\mathcal H}}}

\begin{document}

\title{Hausdorff dimension in a family of self-similar groups}

\author{Zoran \v Suni\'c}
 \thanks{Partially supported by NSF grant DMS-0600975}
 \address{Department of Mathematics, Texas A\&M University, MS-3368, College Station, TX 77845-3368, USA}
 \email{sunic@math.tamu.edu}
 \date{}

\begin{abstract}
For each prime $p$ and a monic polynomial $f$, invertible over $p$,
we define a group $G_{p,f}$ of $p$-adic automorphisms of the $p$-ary
rooted tree. The groups are modeled after the first Grigorchuk
group, which in this setting is the group $G_{2,x^2+x+1}$. We show
that the constructed groups are self-similar, regular branch groups.
This enables us to calculate the Hausdorff dimension of their
closures, providing concrete examples (not using random methods) of
topologically finitely generated closed subgroups of the group of
$p$-adic automorphisms with Hausdorff dimension arbitrarily close to
1. We provide a characterization of finitely constrained groups in
terms of the branching property, and as a corollary conclude that
all defined groups are finitely constrained. In addition, we show
that all infinite, finitely constrained groups of $p$-adic
automorphisms have positive and rational Hausdorff dimension and we
provide a general formula for Hausdorff dimension of finitely
constrained groups. Further ``finiteness'' properties are also
discussed (amenability, torsion and intermediate growth).
\end{abstract}

\maketitle


\section{Introduction}

The group known as the first Grigorchuk group was constructed
in~\cite{grigorchuk:burnside} as a particularly simple example of a
finitely generated infinite torsion group (Burnsde group) and has
since provided solutions to several outstanding problems in
mathematics. For example, this group was the first example of a
group of intermediate
growth~\cite{grigorchuk:growth,grigorchuk:gdegree} and the first
example of amenable but not elementary amenable
group~\cite{grigorchuk:gdegree}. It served as a prime example and
inspiration leading to the definition of many important classes of
groups, such as branch groups, self-similar groups, finitely
constrained groups, groups with $L$-presentations, groups generated
by bounded automata, contracting groups, and so on (for more details
see~\cite{grigorchuk:jibg,bartholdi-g-n:fractal,bartholdi-g-s:branch,sidki:circuit,grigorchuk:unsolved,nekrashevych:b-selfsimilar}).

There have been some generalizations, such as those
in~\cite{grigorchuk:gdegree,grigorchuk:pgps} as well as the so
called spinal groups in~\cite{bartholdi-s:wpg,bartholdi-g-s:branch}.
However, even as the family of groups of Grigorchuk type grew, the
apparent simplicity and symmetry of the first example has not been
repeated (perhaps with the exception of the Gupta-Sidki
example~\cite{gupta-s:burnside}, which actually has its own
interesting properties, and forms a basis for another family of
interesting examples). One of the purposes of this article is to
explicitly construct, for each prime $p$, groups that resemble the
original example as much as possible. In a sense, we are after
siblings of the first Grigorchuk group, not just some distant
relatives. In what follows, one group $G_{p,f}$ will be constructed
for each prime $p$ and each monic polynomial $f$, invertible over
$p$, of degree $m \geq 2$. In this setting, the Grigorchuk group is
the group $G_{2,x^2+x+1}$, defined by the smallest possible prime
$p=2$ and the unique primitive polynomial over $GF(2)$ of the
smallest possible degree $m=2$. We will show that all constructed
examples are self-similar, regular branch groups and we will exhibit
their branching structure (Theorem~\ref{branch}).

Ab\'ert and Vir\'ag showed in~\cite{abert-v:dimension} that the
closure of a group of $p$-adic tree automorphisms of the $p$-ary
rooted tree generated by three randomly chosen automorphisms has
Hausdorff dimension 1 with probability 1. However, no concrete
examples of topologically finitely generated closed groups of
$p$-adic automorphisms of Hausdorff dimension 1 are known. The
branching structure of $G_{p,f}$ enables us to show that $G_{p,f}$,
defined by a polynomial $f$ of degree $m \geq 2$, has Hausdorff
dimension $1-t/p^{m+1}$, where $t=p$, when $p \neq 3$, and $t=3$,
when $p=2$ (Theorem~\ref{dimension}). Thus, for any prime $p$, we
have concrete examples of topologically finitely generated closed
subgroups of the group of $p$-adic automorphisms with Hausdorff
dimension arbitrarily close to 1.

Finitely constrained groups were introduced by Grigorchuk
in~\cite{grigorchuk:unsolved} (under the name ``groups of finite
type''), where it is shown that the closure of the first Grigorchuk
group in the group of automorphisms of the binary tree is an example
of such a group. In the same work (\cite{grigorchuk:unsolved},
Proposition~7.5) Grigorchuk showed that the finitely constrained
groups are always closures of self-similar, regular branch groups,
branching over some level stabilizer. We provide a converse to this
claim, thus characterizing finitely constrained groups in terms of
their branching structure (Theorem~\ref{characterization}). This
immediately implies that the groups constructed here are finitely
constrained groups. Moreover, we show that every infinite, finitely
constrained group of $p$-adic automorphisms has positive and
rational Hausdorff dimension and provide a general formula for the
Hausdorff dimension of finitely constrained groups
(Theorem~\ref{generaldimension}). This provides partial answers to
Problem~7.1.(ii) and Problem~7.1.(iii)
from~\cite{grigorchuk:unsolved}.

At the very end, we discuss some further properties of the
constructed groups, such as amenability, torsion and growth. In
particular, we characterize those examples that are $p$-groups
(Proposition~\ref{torsion}).


\section{Self-similar groups}

We provide a quick and informal introduction to the notion of a
self-similar group.

Let $X$ be a finite alphabet. Our choice of a standard alphabet on
$k$ letters is $X=\{0,1,\dots,k-1\}$. The set $X^*$ of all finite
words over $X$ can be given the structure of a $k$-ary rooted tree.
The empty word $\emptyset$ is the root, the set of words $X^n$ of
length $n$ constitute level $n$ in the tree and each vertex $u$ on
level $n$ has $k$ children, namely $ux$, $x \in X$. The group of
tree automorphisms of the tree $X^*$ is denoted by $\Aut(X^*)$.

The group $\Aut(X^*)$ decomposes algebraically as
\begin{equation}\label{decomposition}
 \Aut(X^*) =
 \Sym_k \ltimes (\Aut(X^*) \times \dots \times \Aut(X^*)),
\end{equation}
where $\Sym_k$, the symmetric group on $X$, acts on $\Aut(X^*)
\times \dots \times \Aut(X^*)= \Aut(X^*)^X$ by permuting the
coordinates. Thus $\Aut(X^*)$ is the permutational wreath product
$\Aut(X^*) = \Sym_k \wr \Aut(X^*)$, i.e., $\Aut(X^*)$ is the
infinitely iterated permutational wreath product $\Aut(X^*) = \Sym_k
\wr \Sym_k \wr \dots$. The normal subgroup $\Aut(X^*) \times \dots
\times \Aut(X^*)$ in the decomposition~\eqref{decomposition} is the
stabilizer of the first level of the tree $X^*$.

The decomposition~\eqref{decomposition} indicates that every tree
automorphisms can be written (decomposed) uniquely as
\begin{equation}\label{gdecomposition}
 g = \pi_g \ (g_0,g_1,\dots,g_{k-1}).
\end{equation}
The decomposition~\eqref{gdecomposition} describes (decomposes!) the
action of $g$ on $X^*$ as follows. The tree automorphisms $g_x$, $x
\in X$, act on the subtrees $xX^*$, $x \in X$, and then the
permutation $\pi_g \in \Sym_k$ permutes these subtrees. The
permutation $\pi_g$ is called the root activity of $g$ and the
automorphisms $g_x$, $x \in X$ are called (the first level) sections
of $g$. The notion of a section can be extended to any vertex,
recursively by $g_\emptyset=g$, $g_{ux} = (g_u)_x$, for $u$ a word
over $X$ and $x$ a letter in $X$. The notion of activity is also
extended to all vertices by declaring that the root activity of the
section $g_u$ is the activity of $g$ at the vertex $u$. The section
of $g$ at $u$ is the unique tree automorphism $g_u$ such that
\[ g(uw) = g(u)g_u(w), \]
for all $w$ in $X^*$. Thus, the section $g_u$ describes the action
of $g$ on the tails of words that start in $u$. In particular, the
decomposition~\eqref{gdecomposition} indicates that the action of
$g$ on the tree $X^*$ is recursively given by
\[ g(xw)  = \pi_g(x)g_x(w), \]
for $x$ a letter in $X$ and $w$ a word over $X$.

\begin{definition}
A group $G$ of tree automorphisms of $X^*$ is self-similar if every
section of every element in $G$ is an element in $G$.
\end{definition}

For a group of tree automorphisms $G \leq \Aut(X^*)$, denote by
$G_n$  the stabilizer of level $n$. Denote by $X_{[n]}$ the finite
subtree of $X^*$ consisting of levels 0 through $n$ in $X^*$. The
quotient $\Aut(X^*)/\Aut(X^*)_n$ is canonically isomorphic to the
group $\Aut(X_{[n]})$ of tree automorphisms of the finite $k$-ary
tree $X_{[n]}$. The group $\Aut(X^*)$, being the inverse limit of
the finite groups $\Aut(X_{[n]})=\underbrace{\Sym_k \wr \Sym_k \wr
\dots \wr \Sym_k}_n$, is a pro-finite group. A natural metric on
$\Aut(X^*)$, derived from the filtration of $\Aut(X^*)$ by its level
stabilizers, is defined by
\[
 d(f,g) =
 \inf \left\{\ \frac{1}{[\Aut(X^*):\Aut(X^*)_n]} \mid
              f^{-1}g \in \Aut(X^*)_n \ \right\}.
\]
Thus, the distance directly depends on the number of levels on which
$f$ and $g$ agree.

Of particular interest is the group of $p$-adic automorphisms, for
$p$ a prime number. Let $\pi = (0 \ 1 \ \dots \ p-1)$ be the
standard cycle on $X=\{0,\dots,p-1\}$. A $p$-adic automorphism $g$
is an automorphism of the $p$-ary rooted tree such that the activity
of $g$ at any vertex is a power of the cyclic permutation $\pi$. The
group of all $p$-adic automorphisms is denoted by $\A(p)$. It a
pro-$p$-group isomorphic to the infinitely iterated wreath product
$C_p \wr C_p \wr \dots$. The group $\A(p)$ is clearly self-similar.
As an example of a $p$-adic automorphism, define
\begin{equation}\label{a}
 a = \pi \ (1,1,\dots,1).
\end{equation}
This automorphism just rigidly permutes all the subtrees on the
first level of $X^*$. i.e., it only affects the first letter in each
word. The group $\A(p)$ of $p$-adic automorphisms is a closed
subgroup of the group of automorphisms of the $p$-ary tree and, as
such, it inherits the metric from $\Aut(X^*)$. However, it is
customary to work with the natural metric derived from the
filtration of $\A(p)$ by its own level stabilizers and defined by
\begin{equation}\label{metric-a}
 d(f,g) =
 \inf \left\{\ \frac{1}{[\A(p):\A(p)_n]} \mid f^{-1}g \in \A(p)_n \ \right\}.
\end{equation}
Note that, when $p=2$, there is no difference between automorphisms
and 2-adic automorphisms, i.e., $\Aut(\{0,1\}^*) = \A(2)$.

Let $G$ be a group of tree automorphisms. For a vertex $u$, the map
$\varphi_u: G_u \to \Aut(X^*)$ mapping each element of the vertex
stabilizer $G_u$ of $u$ in $G$ to its section $g_u$ at $u$ is a
homomorphism. The map $\psi: G_1 \to \Aut(X^*) \times \dots \times
\Aut(X^*)$ defined by $\psi(g) =
(\varphi_0(g),\dots,\varphi_{k-1}(g))$ is also a homomorphism.

A group $G$ of tree automorphisms is self-replicating if
$\varphi(G_u) = G$, for every vertex $u$ in $X^*$.

A group $G$ of tree automorphisms acts spherically transitively on
$X^*$ if it acts transitively on each level of $X^*$.


\section{Examples}

We start with the standard definition of the first Grigorchuk group
as a self-similar group of automorphisms of the binary tree
(equivalently, as a group of 2-adic automorphisms). Then we observe
some basic elements of the construction, which lead us to the
definition of groups associated with maps between finite vector
spaces. A natural faithfulness condition then leads to the
definition of a single group $G_{p,f}$ for each prime $p$ and a
monic polynomial $f$, invertible over $p$.

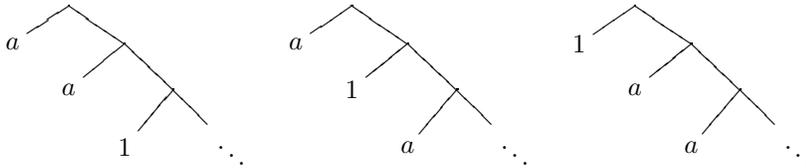
\begin{figure}[!ht]
\[
\xymatrix@C=10pt@R=5pt{
   & *[o][F-]{} \ar@{-}[ld] \ar@{-}[rd] &  & & &
   & *[o][F-]{} \ar@{-}[ld] \ar@{-}[rd] &  & & &
   & *[o][F-]{} \ar@{-}[ld] \ar@{-}[rd] &  & &
 \\
   a &  & *[o][F-]{} \ar@{-}[ld] \ar@{-}[rd] & & &
   a &  & *[o][F-]{} \ar@{-}[ld] \ar@{-}[rd] & & &
   1 &  & *[o][F-]{} \ar@{-}[ld] \ar@{-}[rd] & & &
 \\
   &  a & & *[o][F-]{} \ar@{-}[ld] \ar@{-}[rd] & &
   &  1 & & *[o][F-]{} \ar@{-}[ld] \ar@{-}[rd] & &
   &  a & & *[o][F-]{} \ar@{-}[ld] \ar@{-}[rd] & &
   \\
   &    &  1 & & \ddots &
   &    &  a & & \ddots &
   &    &  a & & \ddots & &
}
\]
\caption{The action of $b$, $c$ and $d$ on the binary
tree}\label{f:bcd}
\end{figure}

The first Grigorchuk group is the self-similar group $G=\langle
a,b,c,d \rangle$ of tree automorphisms generated by the
automorphisms $a$, $b$, $c$, and $d$ defined by
\begin{equation}\label{gg}
 a = (01) \ (1,1), \qquad b = (a,c), \qquad c = (a,d), \qquad d=(1,b).
\end{equation}
The action of $b$, $c$ and $d$ on the binary tree is illustrated
(from left to right) in Figure~\ref{f:bcd}.
The pattern $aa1aa1\dots$ associated to $b$ in Figure~\ref{f:bcd}
extends indefinitely, as do the shifted patterns $a1aa1a\dots$ and
$1aa1aa\dots$ that are associated to $c$ and $d$. Thus, $b$ acts as
$a$ on the left subtree $0X^*$, and it acts as $c$ on the right
subtree $1X^*$. Similarly, $c$ acts as $a$ on the left subtree and
as $d$ on the right one. Finally, $d$ acts trivially on the left
subtree and as $b$ on the right subtree. It is easy to see that $A =
\langle a \rangle \cong \Z/2\Z$ and $B = \langle b,c,d \rangle \cong
\Z/2\Z \times \Z/2\Z$. The cycle $b \stackrel{\varphi_1}{\mapsto} c
\stackrel{\varphi_1}{\mapsto} d \stackrel{\varphi_1}{\mapsto} b$
apparent from the definition of these automorphisms corresponds to
the nontrivial cycle of the automorphism $\rho:B \to B$ defined by
\[
 \rho = \begin{pmatrix} 1 & b & c & d \\ 1 & c & d & b \end{pmatrix}.
\]
Similarly, the definition of the left sections $\varphi_0(b)=a$,
$\varphi_0(c)=a$, $\varphi_0(d)=1$ corresponds to the surjective
homomorphism $\omega:B \to A$ defined by
\[
 \omega = \begin{pmatrix} 1 & b & c & d \\ 1 & a & a & 1 \end{pmatrix}.
\]

Let $p$ be a prime, $A$ and $B$ be the abelian groups $\Z/p\Z$ and
$(\Z/p\Z)^m$, respectively, in multiplicative notation. When
convenient, we also think of $A$ as the field on $p$ elements and of
$B$ as the $m$-dimensional vector space over this field. Let $\rho:B
\to B$ be an automorphism of $B$ (an invertible linear
transformation on the $m$-dimensional vector space $B$) and
$\omega:B \to A$ a surjective homomorphism (nontrivial functional on
$B$).

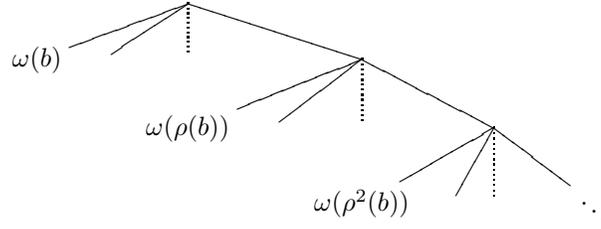
\begin{figure}[!ht]
\[
\xymatrix@C=10pt@R=10pt{
           &  &  *[o][F-]{} \ar@{-}[lld] \ar@{-}[ld] \ar@{..}[d] \ar@{-}[rrd] \\
 \omega(b) &  & & & *[o][F-]{} \ar@{-}[lld] \ar@{-}[ld] \ar@{..}[d] \ar@{-}[rrd]  \\
           &  & \omega(\rho(b)) & & &  &  *[o][F-]{} \ar@{-}[lld] \ar@{-}[ld] \ar@{..}[d] \ar@{-}[rrd]  \\
           &  &                 & & \omega(\rho^2(b)) && &&
           \ddots
}
\]
\caption{The action of $b$ on the tree $X^*$}\label{f:actionb}
\end{figure}

Let $A = \langle a \rangle$ and let $a$ act faithfully on the
$p$-ary rooted tree as the $p$-adic automorphism defined
in~\eqref{a}. Let the action of $b \in B$ on $X^*$ be recursively
defined by
\[ b = (\omega(b),1,\dots,1,\rho(b)). \]
The action of $b$ is illustrated in Figure~\ref{f:actionb}.

\begin{definition}
For each pair $(\omega,\rho)$ of a surjective homomorphism $\omega:B
\to A$ and an automorphism $\rho:B \to B$, define a group
$G_{\omega,\rho}$ of $p$-adic tree automorphisms by
\[ G_{\omega,\rho} = \langle A \cup B \rangle \leq \A(p). \]
\end{definition}

Since each section of a generator of $G_{\omega,\rho}$ is another
generator or is trivial, the defined groups are self-similar. We
will not emphasize this approach, but it is also clear that each
group $G_{\omega,\rho}$ can be defined by a finite automaton
(see~\cite{grigorchuk-n-s:automata} for definition).

We now observe that different pairs $(\omega, \rho)$ may define the
same group.

\begin{proposition}
Let $\alpha:A \to A$ and $\beta:B \to B$ be automorphisms such that
the following diagram commutes:
\[
\xymatrix{
 A \ar@{->}[d]^{\alpha} & B \ar@{->}[l]_{\omega_1} \ar@{->}[r]^{\rho_1} \ar@{->}[d]^{\beta} & B \ar@{->}[d]^{\beta} \\
 A & B \ar@{->}[l]^{\omega_2} \ar@{->}[r]_{\rho_2} & B
}
\]
Then $G_{\omega_1,\rho_1}=G_{\omega_2,\rho_2}$ as subgroups of
$\A(p)$.
\end{proposition}
\begin{proof}
The automorphism $\alpha$ is just a power automorphism given by
$\alpha(x)=x^k$ for some fixed $k=1,\dots,p-1$ (in vector space
terminology, $\alpha$ is multiplication by the non-zero scalar $k$).
Define $\alpha':B \to B$ to be the corresponding power automorphism
of $B$, given by $\alpha'(x)=x^k$. Then
$\alpha\omega_1=\omega_1\alpha'$ and $\alpha'\rho_1=\rho_1\alpha'$
and, for $b \in B$, the automorphisms $\beta(b)$ in
$G_{\omega_2,\rho_2}$ and $\alpha'(b)$ in $G_{\omega_1,\rho_1}$ are
the same automorphism of $X^*$. The latter is clear from the fact
that
\[ \omega_2\rho_2^s(\beta(b)) = \omega_2\beta\rho_1^s(b) =
 \alpha\omega_1\rho_1^s(b) = \omega_1\rho_1^s(\alpha'(b)), \]
for $s \geq 0$.
\end{proof}

Thus in concrete examples we may always choose a conjugate
representative of $\rho$. For example, we may replace a given group
$G_{\rho_2,\omega_2}$ by $G_{\rho_1,\omega_1}$, where $\rho_1$ is
any conjugate $\rho_1=\beta^{-1}\rho_2\beta$ and
$\omega_1=\omega_2\beta$ ($\alpha$ may be chosen to be the identity
map here).

We will provide a simple condition providing faithfulness of the
action of $B$ on $X^*$. We first introduce some notation. For $i \in
\Z$, let $B_i$ be the subgroup of $B$ obtained as the image of the
kernel $\Ker(\omega)$ under the $i$-th iteration of $\rho$,
i.e.~$B_i=\rho^i(\Ker(\omega))$. In particular,
\[ B_0 = \Ker(\omega), \qquad B_1 = \rho(B_0), \qquad B_{-1}=\rho^{-1}(B_0).\]

\begin{proposition}\label{faithful}
The following conditions are equivalent:

(i) the action of $B$ on $X^*$ is faithful.

(ii) no nontrivial orbit of $\rho$ is contained in the kernel of
$\omega$.

(iii) no nontrivial $\rho$-invariant subspace of $B$ is contained in
the kernel of $\omega$.

(iv) $B$ is $\rho$-cyclic, the minimal polynomial $f$ of $\rho$ has
degree $m$, $f(x)=x^m+a_{m-1}x^{m-1} +  \dots + a_0$, and there is a
basis of $B$ with respect to which the matrix of $\rho$ is given by
\begin{equation}\label{mro}
 M_\rho =
 \begin{bmatrix}
 0 &  0 & \dots & 0 & -a_0 \\
 1 &  0 & \dots & 0 & -a_1 \\
 0 &  1 & \dots & 0 & -a_3 \\
 \vdots & \vdots & \ddots & \vdots & \vdots \\
 0 & 0  & \dots & 1 & -a_{m-1},
 \end{bmatrix}
\end{equation}
and the matrix of $\omega$ is given by
\begin{equation}\label{momega}
 M_\omega =
 \begin{bmatrix}
 0 &  0 & \dots & 0 & 1
\end{bmatrix}.
\end{equation}
\end{proposition}

\begin{proof}

(i) is equivalent to (ii). The $\rho$-orbit of $b \in B$ fully
describes the action of $b$ on $X^*$ (see Figure~\ref{f:actionb}).
The action of $b$, $b \neq 1$, is nontrivial if and only if the
$\rho$-orbit of $b$, is not entirely in $B_0$.

The statement (ii), (iii), and (iv) are all trivially fulfilled when
$m=1$. We assume $m \geq 2$.

(ii) implies (iii). If a nontrivial $\rho$-invariant subspace $W$ of
$B$ is contained in $B_0$, then the $\rho$-orbit of any nontrivial
element of $W$ is included in $B_0$.

(iii) implies (iv). Assume that nontrivial $\rho$-invariant subspace
of $B$ is contained in $B_0$. Since the subspaces $B_i$, $i \in \Z$,
have codimension 1 in $B$, which is $m$-dimensional (recall that $m
\geq 2$),
\[ W=B_0 \cap B_1 \cap \dots \cap B_{m-2} \]
has dimension at least 1, and is thus nontrivial. Thus, there exists
a piece
\[ d \mapsto \rho(d) \mapsto \rho^2(d) \mapsto \dots \mapsto \rho^{m-2}(d)\]
of a nontrivial $\rho$-orbit that lies entirely in $B_0$. If we
assume that the minimal polynomial $f$ of $\rho$ has degree strictly
smaller than $m$, then the dimension of the $\rho$-cyclic space
$\langle d \rangle_\rho$ generated by $d$ is at most $m-1$, and
therefore the above piece contains a basis for $\langle d
\rangle_\rho$. This is a contradiction, since this would imply that
the whole $\rho$-cyclic space $\langle d \rangle_\rho$ lies in
$B_0$. Thus the minimal polynomial of $\rho$ has degree $m$,
implying that $B$ is $\rho$-cyclic. Moreover, $B = \langle d
\rangle_\rho$ and a cyclic basis of $B$ consists of
$d,\rho(d),\dots,\rho^{m-2}(d),\rho^{m-1}(d)$. The matrix of $\rho$
with respect to this basis is exactly the matrix in~\eqref{mro},
which is the companion matrix of the minimal polynomial $f(x) = x^m
+ a_{m-1}x^{m-1} + \dots + a_0$ of $\rho$. Moreover, since
$d,\rho(d),\dots,\rho^{m-2}(d) \in W \subseteq B_0$, their images
under $\omega$ are trivial. On the other hand, the image of
$\rho^{m-1}$ must not be trivial, since $\omega$ is not trivial.
Thus the matrix of $\omega$ is the one given in \eqref{momega}
(after re-scaling of the whole cyclic basis, if necessary).

(iv) implies (ii). In our chosen basis, the functional $\omega$ just
reads the last coordinate of an input vector. The only way a
$\rho$-orbit would lie entirely in $B_0$ is if the last coordinate
stayed equal to 0 in the whole orbit. But
\[ M_\rho [x_0 \ x_1 \ \dots \ x_{m-2} \ 0]^T = [0 \ x_0 \ x_1 \ \dots \ x_{m-2}]^T, \]
which forces $x_{m-2}$ to be 0, and then $M_\rho[0 \ x_0 \ \dots \
x_{m-3} \ 0]^T = [0 \ 0 \ x_0 \ \dots \ x_{m-3}]^T$ forces $x_{m-3}$
to be 0, and so on. Eventually we obtain that the zero vector is in
the orbit, and therefore the orbit is trivial. Thus $B_0$ does not
contain nontrivial $\rho$-orbits.
\end{proof}

From now on, only pairs $(\omega,\rho)$ for which the above
faithfulness condition is fulfilled are considered.

Further, we often choose to work with the basis guaranteed by and
described in Proposition~\ref{faithful}. In order to specify a group
in the family, we only specify the prime $p$ and the polynomial
$f=x^m+a_{m-1}x^{m-1}+ \dots + a_1x+a_0$. In particular, there is a
sequence of elements $d_0,d_1,\dots,d_{m-2}$ in $B_0$ and an element
$d_{m-1}$ in $B_1-B_0$, such that
\begin{equation}\label{d0b}
\begin{array}{ccccccl}
 d_0 & = &(1, &1, &\dots, &1, &d_1) \\
 d_1 & = &(1, &1, &\dots, &1, &d_2) \\
 \dots \\
 d_{m-2} & = &(1, &1, &\dots, &1, &d_{m-1}) \\
 d_{m-1} & = &(a,&1,&\dots,&1,&b')
\end{array}
\end{equation}
in $\Aut(X^*)$, where $b' =\rho(d_{m-1}) = d_0^{-a_0} d_1^{-a_1}
\dots d_{m-2}^{-a_{m-2}} d_{m-1}^{-a_{m-1}}$. The action of $B$, and
therefore the whole group $G_{\omega,\rho}=G_{p,f}$, can be
recovered from the relations in~\eqref{d0b}.

\begin{definition}
Let $p$ be a prime and $f=x^m+a_{m-1}x^{m-1}+\dots+a_1x+a_0$ be a
polynomial, which is invertible over $p$ (i.e.~$a_0 \neq 0$). The
group $G_{p,f}$ is the group $G_{\omega,\rho}$, where $\omega$ and
$\rho$ are given by \eqref{momega} and ~\eqref{mro}.
\end{definition}

In practice, the decomposition formulas~\eqref{d0b} are extended to
include the full orbit of $d_0$. We call the elements in the
$\rho$-orbit of $d_0$ standard generators of $B$ and the elements
$d_0,\dots,d_{m-2},d_{m-1}$ standard basis elements of $B$. For
future reference, note that the standard basis elements
$d_0,\dots,d_{m-2}$ generate the kernel $B_0$, while the standard
basis elements $d_1,\dots,d_{m-2},d_{m-1}$ generate its image $B_1$.

We list now all examples that have been explicitly mentioned in the
literature and that fit in our construction.

\begin{example}[The infinite dihedral group] The infinite dihedral
group $D_\infty$ is the group $G_{2,x-1}$. We have
$M_\rho=M_\omega=[1]$ and the action of the only nontrivial element
$b$ of $B$ on $X^*$ is given by
\[ b = (a,b). \]
\end{example}

\begin{example}[The first Grigorchuk group] The first Grigorchuk group
$G$ is the group $G_{2,x^2+x+1}$. The standard matrices of $\rho$
and $\omega$ are
\begin{equation}\label{fibonacci}
  M_\rho = \begin{bmatrix} 0 & 1 \\ 1 & 1 \end{bmatrix}
\end{equation}
and
\[ M_\omega = \begin{bmatrix} 0 & 1 \end{bmatrix}. \]
Setting $d=[1 \ 0]^T$, $b=[0 \ 1]^T$, leads to decomposition
formulas
\[ d = (1,b), \qquad b = (a,db), \]
which are sufficient to define the action on $X^*$. However, in
practice, the decomposition formulas are given for the whole
$\rho$-orbit $d \mapsto b \mapsto c \mapsto d$
\[ d = (1,b), \qquad b = (a,c), \qquad c = (a,d), \]
where $c=db=[1 \ 1]^T$, which is the usual way used to describe the
action of Grigorchuk group (see~\eqref{gg}). The pattern
$aa1aa1\dots$ and its shifts, evident in the action of $b$, $c$ and
$d$ on $X^*$ (see Figure~\ref{f:bcd}), is just a manifestation of
the modulo 2 Fibonacci sequence $110110\dots$, which can be defined
by iterations of the Fibonacci matrix $M_\rho$ in~\eqref{fibonacci}.
\end{example}

\begin{example}[Grigorchuk overgroup] The Grigorchuk overgroup is
the group $G_{2,x^3+1}$. If the standard matrices
\[ M_\rho = \begin{bmatrix}
 0 & 0 & 1 \\
 1 & 0 & 0 \\
 0 & 1 & 0 \end{bmatrix} \]
and
\[ M_\omega = \begin{bmatrix} 0 & 0 & 1 \end{bmatrix} \]
are conjugated appropriately (i.e., we introduce change of basis by
$d_0 \mapsto d_0d_1d_2$, $d_1 \mapsto d_0d_1$, $d_2 \mapsto d_1d_2$)
we get
\[ M_\rho' = \begin{bmatrix}
 1 & 0 & 0 \\
 0 & 0 & 1 \\
 0 & 1 & 1 \end{bmatrix} \]
 and
\[ M_\omega' = \begin{bmatrix} 1 & 0 & 1 \end{bmatrix}, \]
which clearly shows that $G_{2,x^3+1}$ is generated by a copy of the
dihedral group (generated by $a$ and $d_0d_1d_2$) and a copy of the
first Grigorchuk group (generated by $a$, $d_0d_1$ and $d_1d_2$).
\end{example}

\begin{example}[Grigorchuk-Erschler group]
Grigorchuk-Erschler group was first defined
in~\cite{grigorchuk:gdegree}, as a member of an uncountable family
of groups of intermediate growth. This group is studied in more
detail in~\cite{erschler:subexponential}. It is the group
$G_{2,x^2+1}$.
\end{example}

\begin{example}[The Fabrykowski-Gupta group] The Fabrykowski-Gupta
group was introduced and studied in~\cite{fabrykovski-g:growth}. It
is the group $G_{3,x-1}$. More generally, the group $G_{p,x-1}$ has
been studied in~\cite{grigorchuk:jibg}.
\end{example}

\begin{proposition}
Let $f$ be a monic polynomial, invertible over $p$, which factors as
$f=f_1f_2$, for some nonconstant monic polynomials. Then

(a) $G_{p,f_1} \leq G_{p,f}$.

(b) $G_{p,f}$ is generated by its subgroups $G_{p,f_1}$ and
$G_{p,f_2}$, provided $f_1$ and $f_2$ are relatively prime.
\end{proposition}
\begin{proof}
(a) Consider the element $b=f_2(\rho)(d_0)$. The minimal polynomial
of the $\rho$-cyclic subspace $W_b=\langle b \rangle_\rho$ is equal
to $f_1$. Moreover, if the degree of $f_1$ is $m_1$, the elements
$b,\rho(b),\dots,\rho^{m_1-2}(b)$ are in the kernel of $\omega$,
while $\rho^{m_1-1}(b)$ is not. Thus, the subgroup of
$G_{p,f}=\langle A \cup B \rangle$ generated by $A$ and $W_b$ is the
group $G_{p,f_1}$.

(b) If $f=f_1f_2$, where $f_1$ and $f_2$ are relatively prime then
the $\rho$-cyclic basis of the $\rho$-cyclic subspace $W_{b_1}$ with
minimal polynomial $f_1$ and the $\rho$-cyclic basis of the
$\rho$-cyclic subspace $W_{b_2}$ with minimal polynomial $f_2$ form
a basis of $B$. Thus $G$ is generated by $A$, $W_{b_1}$ and
$W_{b_2}$, i.e., $G$ is generated by $G_{p,f_1}$ and $G_{p,f_2}$.
\end{proof}

We will see that, for a fixed $p$ and $m$, all groups $G_{p,f}$
defined by a polynomial of degree $m$ share many properties. For a
fixed prime $p$, we denote the family of groups defined by
polynomials of degree $m$ over the field of order $p$ by $\G_{p,m}$.
We denote by $\G$ the union of the families $\G_{p,m}$ for all
primes and all degrees.

We note again that the groups of the form $G_{p, x-1}$ were studied
in~\cite{grigorchuk:jibg}. From now on, we mostly restrict our
attention to the case $m \geq 2$, which is in many ways crucially
different than the case $m=1$.


\section{Elementary properties}

\begin{proposition}\label{commutator}
Let $G$ be a group in $\G$. The maps $\pi_A:G \to A$ and $\pi_B: G
\to B$ given by
\[ \pi_A(a_1b_1 \cdots a_kb_ka_{k+1}) = a_1 \cdots a_{k+1} \]
and
\[ \pi_B(a_1b_1 \cdots a_kb_ka_{k+1}) = b_1 \cdots b_k, \]
for $a_1, \dots, a_{k+1} \in A$ and $b_1, \dots, b_k \in B$, are
well defined homomorphisms. The kernel of $\pi_A$ is the first level
stabilizer $G_1$, and it is equal to the normal closure $B^G$ of $B$
in $G$. The kernel of $\pi_B$ is the normal closure $A^G$ of $A$ in
$G$.

The commutator of $G$ is
\[ G' = \Ker(\pi_A) \cap \Ker(\pi_B) \]
and the abelianization $G/G'$ is isomorphic to $A \times B$.
\end{proposition}
\begin{proof}
We provide only a sketchy proof, since this statement is analogous
to the corresponding statements about the spinal groups defined in
~\cite{bartholdi-s:wpg}.

The claims related to $\pi_A$ follow from the fact that $B \leq G_1$
and that the element $g=a_1b_1 \cdots a_kb_ka_{k+1}$ can be
rewritten as
\begin{equation}\label{g=ab}
 g = a_1 \dots a_{k+1} b_1^{a_2\dots a_{k+1}} \dots b_{k-1}^{a_ka_{k+1}}b_k^{a_{k+1}}.
\end{equation}

To prove the claims related to $\pi_B$ we need to show that if
$[a_1]b_1 \cdots a_kb_k[a_{k+1}]$ is a word over $S=A \cup B
\setminus \{1\}$ (the bracketed letters in the front and at the end
may be absent) that represents the identity in $G$, then $b_1 \dots
b_k=1$ in $B$. This follows by induction on the length $n$ of the
words in $S$. The base cases regarding length 0 and 1 are trivial.
For the inductive step, first observe that
\begin{equation} \label{contraction}
 |g_i| \leq \frac{|g|+1}{2},
\end{equation}
for $i=0,\dots,p-1$, and that each letter $b_j$, $j=1,\dots,k$
contributes $\rho(b_j)$ to exactly one of the sections $g_i$,
$i=0,\dots,p-1$. Thus if $b_1 \dots b_k \neq 1$ in $B$ then
$\rho(b_1)\dots\rho(b_k) \neq 1$ in $B$ and therefore at least one
of the words over $S$ representing the sections of $g$ does not
project to 1 in $B$. But this contradicts the inductive assumption,
since the words representing the sections must represent 1 in $G$
and are shorter than $n$, as long as $n \geq 2$
(by~\eqref{contraction}).

Since the images of both $\pi_A$ and $\pi_B$ are abelian groups, we
must have $G' \leq \Ker(\pi_A) \cap \Ker(\pi_B)$. On the other hand,
any element $g=[a_1]b_1 \cdots a_kb_k[a_{k+1}]$ in $G$ can be
rewritten as  $\pi_A(g)\pi_B(g)h$, where $h \in G'$. Thus if $g \in
\Ker(\pi_A) \cap \Ker(\pi_B)$, then $g \in G'$.
\end{proof}

\begin{proposition}\label{subdirect}
Let $G$ be a group in $\G$.

(a) The map $\psi$ is a subdirect embedding of $G_1$ into $G \times
\dots \times G$

(b) $G$ is a self-replicating, spherically transitive group.
\end{proposition}
\begin{proof}
(a) Since, for $b \in B$,
\[ \rho^{-1}(b) = (\omega(\rho^{-1}(b)),1,\dots,1, b), \]
and, for the standard basis generator  $d_{m-1}$,
\[ d_{m-1}^a = (1,\dots,1,b',a) \]
we see that the the coordinate $p-1$ in $\psi(G_1)$ projects onto
$G$. Conjugating by $a$, we can move any coordinate in $\psi(G_1)$
to any position. Thus $\psi: G_1 \to G \times \dots \times G$ is a
subdirect embedding (surjective on each coordinate).

(b) For a vertex $u=x_1 \dots x_n$, the map $\varphi_u: G_u \to G$
is equal to the composition $\varphi_{x_n} \dots \varphi_{x_1}$ and
each of the maps $\varphi_x$, $x \in X$ is surjective, the map
$\varphi_u$ is surjective as well. Thus $G$ is self-replicating.

A self-replicating group that acts transitively on level 1 is always
spherically transitive.
\end{proof}


\section{Branching property}

Some notable examples of finitely generated branch groups, such as
the first Grigorchuk group~\cite{grigorchuk:burnside} and
Gupta-Sidki group~\cite{gupta-s:burnside}, were known and studied
for two decades before some of their crucial properties were
isolated and served as the model for the notion of a branch group,
which was formally introduced in~\cite{grigorchuk:jibg}. The
relatively late arrival of the formal definition came despite the
fact that the branching property was known and used for both
examples for a long time. The branching property was established for
the first Grigorchuk group in~\cite{grigorchuk:gdegree}, while for
the Gupta-Sidki group this was done in~\cite{gupta-s:extension}.
According to our needs, we will only define the notion of a regular
branch group (see~\cite{grigorchuk:jibg,bartholdi-g-s:branch} for
more information).

Let $H^{(0)},\dots,H^{(k-1)}$ and $H$ be groups of automorphisms of
the $k$-ary rooted tree. We say that $H$ geometrically contains
$H^{(0)} \times \dots \times H^{(k-1)}$ if $H^{(0)} \times \dots
\times H^{(k-1)} \leq \psi(H_1)$. Thus, for each $x \in X$, $H$
contains a subgroup that acts on the subtree $xX^*$ exactly as the
group $H^{(x)}$ does on $X^*$ and acts trivially on all other
subtrees on the first level.

\begin{definition}
A group $G$ of tree automorphisms is a regular branch group,
branching over its subgroup $K$ if

(i) $G$ acts spherically transitively on $X^*$

(ii) $K$ is a normal subgroup of finite index in $G$

(iii) $K$ geometrically contains $K \times \dots \times K$.
\end{definition}

We will show that our examples are regular branch groups, with the
exception of the dihedral group $D_\infty=G_{2,x-1}$
(Lemma~\ref{overg'}, for $p \neq 3$, and Lemma~\ref{overk}, for
$p=2$).

\begin{theorem}\label{branch}
(a) Let $G$ be a group in $\G_{p,m}$, $p \neq 2$, $m \geq 2$. Then
$G$ is a regular branch group over the commutator subgroup $G'$.

The index of $G'$ in $G$ is $[G:G']=p^{m+1}$, while the index of
$\psi^{-1}(G' \times \dots \times G')$ in $G'$ is $p^{m(p-1)}$.

(b) Let $G$ be a group in $\G_{2,m}$, $m \geq 2$. Then $G$ is a
regular branch group over
\[ K = \langle [a,b] \mid b \in B_1 \rangle^G. \]

The index of $K$ in $G$ is $[G:K]=2^{m+2}$, while the index of
$\psi^{-1}(K \times \dots \times K)$ in $K$ is $2^{m}$.
\end{theorem}

\subsection{Odd prime case}

For the duration of this subsection, containing a proof of
Theorem~\ref{branch}(a), we assume that $G$ is a group in
$\G_{p,m}$, for some $p \neq 2$, $m \geq 2$. The conditions that $p$
is an odd prime or that $m \geq 2$ are not needed all the time, so
we clearly indicate each time which conditions are used. This is
done with the intention of recycling some of the results and using
them in the case $p=2$.

\begin{lemma}\label{overg'}
Let $G$ be a group in $\G_{p,m}$, $p \neq 2$. Then $G$ is a regular
branch group over the commutator subgroup $G'$.
\end{lemma}
\begin{proof}
Let $b$ be an element in $B$ with $\omega(b)=a$. For an arbitrary
element $b' \in B$, $[b^a,b']$ is in $G'$ and
\[
 \begin{array}{rrrccccccl}
 \psi( [(b^a,b'] ) & = &( &1, &1, & \dots, & 1, & \rho(b^{-1}), &\omega(b^{-1}) &)\\
                   &   &( &\omega((b')^{-1}), &1, & \dots,&1,&1,&\rho((b')^{-1})&)\\
                   &   &( &1,&1, & \dots, & 1, & \rho(b), &\omega(b)&)\\
                   &   &( &\omega(b'),&1,&\dots,&1,&1,&\rho(b')&) =\\
                   & = &( &1, &1, &\dots,&1,&1, &[a,\rho(b')]&).
\end{array}
\]
Note that the condition $p \neq 2$ is used in the above calculation.
Thus, for $b' \in B$,
\[ \psi( [(b^a,\rho^{-1}(b')] ) = (1,\dots,1,1,[a,b']) \in \psi(G'). \]
Since $\psi$ is a subdirect embedding (Proposition~\ref{subdirect}),
for every $g \in G$, there exists $h \in G_1$ with
$\varphi_{p-1}(h)=g$ and therefore, for every $g \in G$,
\[ \psi([b^a,\rho^{-1}(b')]^h)=(1,\dots,1,[a,b']^g) \in \psi(G'). \]
This implies that $1 \times \dots 1 \times G'\leq \psi(G')$ and,
since $G$ acts transitively on level 1, $G' \times \dots \times G'
\leq \psi(G')$.

Since $G'$ is normal, has finite index in $G$
(Proposition~\ref{commutator}) and $G$ acts spherically transitively
(Proposition~\ref{subdirect}) it follows that $G$ is a branch group
over $G'$.
\end{proof}

We emphasize that the above proof is valid for $m=1$, which will not
be the case with most of our other proofs that essentially use the
existence of a nontrivial kernel $B_0$ of $\omega$.

\begin{lemma}
Let $G$ be a group in $\G_{p,m}$, $m \geq 2$. For every nontrivial
element $d \in B_0$, the subgroup $D=\langle a,d \rangle \leq G$ is
finite. Moreover $D$, is isomorphic to the the wreath product
$\Z/p\Z \wr \Z/p\Z$ (and therefore has $p^{p+1}$ elements).
\end{lemma}
\begin{proof}
Both $a$ and $d$ have order $p$. Further, for $i=0,\dots,p-1$,
\[ d^{a^{i}} = (1,1,\dots,\rho(d),1,\dots,1), \]
where the section $\rho(d)$ appears at the coordinate $p-1-i$. Thus,
$d^{a^{i}}$, $i=0,\dots,p-1$, commute and the normal closure
$\langle d \rangle^D$ of $d$ in $D$ is the elementary $p$ group
$(\Z/p\Z)^p$ on which $a$ acts by shifting coordinates. The
intersection $\langle a \rangle \cap \langle d \rangle^D $ is
trivial since $\langle d \rangle^D$ stabilizes $X$ (level 1 in the
tree), while $A$ acts freely on $X$.
\end{proof}

\begin{lemma}\label{lemma:DB_1}
Let $G$ be a group in $\G_{p,m}$, $m \geq 2$. For any element $d$ in
$B_0 - B_1$, the group $G$ decomposes as
\[ G = D \ltimes \overline{B_1}, \]
where $D = \langle a,d \rangle$ and $\overline{B_1} = B_1^G$ is
the normal closure of $B_1$ in $G$.
\end{lemma}
\begin{proof}
First, note that $B_0 - B_1$ is not empty (the standard generator
$d_0$ is in this set). The group $G$ is generated by $d$, $B_1$, and
$a$. Thus $G = D\overline{B_1}$. It remains to be shown that $D \cap
\overline{B_1} = 1$. Assume that $g \in D \cap \overline{B_1}$.
Since $g$ is an element of $D$ that stabilizes the first level, it
has the form $g=(d_1^{n_0},\dots,d_1^{n_{p-1}})$, where $d_1 =
\rho(d)$. Let $x \in X$. Since $g \in \overline{B_1}$, the section
$g_x$ of $g$ has the form $(a_1^{h_1}\rho(b_1)^{h_1'} \dots
a_k^{h_k} \rho(b_k)^{h_k'} a_{k+1}^{h_{k+1}})$, where $a_i \in A$,
$b_i \in B_1$, $h_i \in G$, and $h_i' \in G$ for all appropriate
indices. This implies that $\pi_B(g_x) = \rho(b_1) \dots \rho(b_k) =
\rho(b_1 \dots b_k) \in \rho(B_1) = B_2$. On the other hand, since
$g_x = d_1^{n_x}$, we obtain that $d_1^{n_x} \in B_2$. However, $d
\in B_0-B_1$, and therefore $d_1 \in B_1 - B_2$. Thus, $d_1 \not\in
B_2$ and, since the order of $d_1$ is $p$, we must have $n_x=0$.
Since each section of $g$ must be trivial, we obtain that $g=1$.
\end{proof}

In order to determine the index of $\psi^{-1}(G' \times \dots \times
G')$ in $G'$ we need to determine the index of $\psi(G_1)$ in $G
\times \dots \times G$.

\begin{lemma}
Let $G$ be a group in $\G_{p,m}$, $p \neq 2$, $m \geq 2$. Then there
exists an element $c \in B_{-1} - B_0$ and an element $d \in B_0 -
B_1$ such that
\[ c = (a,1,\dots,1,d). \]
Let
\[\hat{C} = \langle c, c^a, \dots, c^{a^{p-1}} \rangle =  \langle (a,1,\dots,1,d), (1,\dots,1,d,a), \dots,
 (d,a,1,\dots,1) \rangle. \]
Then,

(a)
\[ \psi(G_1) =
 \hat{C} \ltimes (\underbrace{\overline{B_1} \times \cdots \times \overline{B_1}}_p).
\]

(b)
\[ \underbrace{G \times \cdots \times G}_p =
  (\hat{D} \ltimes \hat{C}) \ltimes (\underbrace{\overline{B_1} \times \cdots \times \overline{B_1}}_p),
\]
where
\[ \hat{D} = \langle (d,1,\dots,1), \dots, (1,\dots,1,d) \rangle. \]

(c) the index of $\psi(G_1)$ in $G \times \dots \times G$ is $p^p$.

(d) the index of $\psi^{-1}(G' \times \dots \times G')$ in $G'$ is
$p^{m(p-1)}$.
\end{lemma}
\begin{proof}
Let $d_0$ be the first standard generator and $c' = \rho^{-1}(d_0)
\in B_{-1} - B_0$. Let $c'=(a^k,1,\dots,1,d_0)$, for some $k
\in\{1,\dots,p-1\}$ ($k$ cannot be 0, since $c' \not\in B_0$). Let
$\ell$ be the multiplicative inverse of $k$ modulo $p$, $c =
(c')^\ell$ and $d=\rho(c)=d_0^\ell$. Then $c = (a,1,\dots,1,d)$, $c
\in B_{-1} - B_0$ and $d \in B_0 - B_1$.

(a) The fact that $\hat{C} \cap (\overline{B_1} \times \dots \times
\overline{B_1}) = 1$ follows from the fact that $D \cap
\overline{B_1} = 1$ (Lemma~\ref{lemma:DB_1}) and $\hat{C} \leq D
\times \dots \times D$ .

The stabilizer $G_1$, which is the closure of $B$ in $G$, is
generated by
\[ \{b^{a^k}|b \in B_0,k=0,\dots,p-1\} \cup
 \{c^{a^k}|k=0,\dots,p-1\}. \]
Since the group $\hat{C} \ltimes (\overline{B_1} \times \dots \times
\overline{B_1})$ contains all images under $\psi$ of these
generators we see that $\psi(G_1) \leq \hat{C} \ltimes
(\overline{B_1} \times \dots \times \overline{B_1})$. For arbitrary
$g \in G$ there exists $h \in G_1$ with $\varphi_{p-1}(h)=g$. Thus,
for $b \in B_1$, $\psi(\rho^{-1}(b)^h) = (1,\dots,1,b^g) \in
\psi(G_1)$, which shows that $1 \times \cdots \times 1 \times
\overline{B_1} \leq \psi(G_1)$. Similarly, for $i=0,\dots, p-1$,
$\psi(\rho^{-1}(b)^{ha^i}) = (1,\dots,1,b^g,1\dots,1)$, where the
nontrivial coordinate appears at $p-1-i$. Thus $\overline{B_1}
\times \dots \times \overline{B_1} \leq \psi(G_1)$. Since $\hat{C}
\leq \psi(G_1)$, it follows that $\hat{C} \ltimes (\overline{B_1}
\times \dots \times \overline{B_1}) \leq \psi(G_1)$.

(b) Since $G = D \ltimes \overline{B_1}$, all we need to show is
that $D \times \dots \times D = \hat{D} \ltimes \hat{C}$. We have
(note the importance of the assumption $p \neq 2$),
\begin{multline*}
 (1,\dots,1,[a,d]) = \\ = (1,\dots,1,d^{-1},a^{-1})(a^{-1},1,\dots,1,d^{-1})
 (1,\dots,d,a)(a,1,\dots,1,d) = \\
 = [c^a,c] \in \hat{C}.
\end{multline*}
For any $g$ in $D$ there exists $h$ in $\hat{C}$ such that
$\varphi_{p-1}(h)=g$. Thus, for any $g \in D$, $\hat{C} \ni
[c^a,c]^h = (1,\dots,1,[a,d]^g)$, which shows that $1 \times \dots
\times 1 \times D' \leq \hat{C}$. Since $\hat{C}$ is closed under
conjugation by $a$ it follows that $D' \times \dots \times D' \leq
\hat{C}$.

The subgroup $\hat{C} \leq D \times \dots \times D$ can be described
as follows. Let $\tilde{C}$ be the set of elements
$g=(g_0,\dots,g_{p-1})$ in $D \times \dots \times D$ such that, for
$i=0,\dots,p-1$, the exponent of $d$ in $\pi_B(g_i)$ is equal to the
exponent of $a$ in $\pi_A(g_{i+1})$ (indices modulo $p$). We claim
that $\hat{C}=\tilde{C}$. Indeed, $\tilde{C}$ is closed for products
and contains all generators of $\hat{C}$. Thus $\hat{C} \leq
\tilde{C}$. On the other hand, an arbitrary element $g$ in
$\tilde{C}$ can be written as
\[
 g =
(a^{n_{p-1}}d^{n_0},d^{n_1}a^{n_0},d^{n_2}a^{n_1},\dots,d^{n_{p-1}}a^{n_{p-2}})h,
\]
where $h \in D' \times \dots \times D' \leq \hat{C}$. Since
$(a^{n_{p-1}}d^{n_0},d^{n_1}a^{n_0},\dots,d^{n_{p-1}}a^{n_{p-2}}) =
c^{n_{p-1}} (c^a)^{n_{p-2}} \dots (c^{a^{p-1}})^{n_0}$, we see that
$\tilde{C} \leq \hat{C}$.

It is clear now that $D \times \dots \times D = \hat{D}\hat{C}$,
since an arbitrary element in $D \times \dots \times D$ can be
multiplied by an element in $\hat{D}$ in such a way that the
obtained product satisfies the requirements describing $\tilde{C}$.
Also, it is clear that $\hat{D} \cap \hat{C} =1$. Indeed,
$(d^{n_0},\dots,d^{n_{p-1}}) \in \hat{C}=\tilde{C}$ implies
$n_0=\dots=n_{p-1}=0$.

Conjugating a coordinate of an element $g$ in $\tilde{C}$ by $a$ or
by $d$ does not affect the condition describing $\tilde{C}$. Thus
$\tilde{C}=\hat{C}$ is a normal subgroup of $D \times \dots \times
D$. Therefore $D \times \dots \times D = \hat{D} \ltimes \hat{C}$.

(c) Follows from (a) and (b) and the observation that $|\hat{D}| =
p^p$.

(d) All the relevant subgroups and indices are given in
Figure~\ref{f:indices}.
\begin{figure}[!ht]
\[
\xymatrix{
 \ar@{<..>}[dd]^{p^{m+1}} &&
 G  & G \times \dots \times G && \ar@{<..>}[ddd]^{p^{pm+p}}\\
 && G_1  \ar@{-}[u]^p \ar@{->}[r]^{\psi}         & \psi(G_1) \ar@{-}[u]_{p^p} && \\
 && G'  \ar@{-}[u]^{p^{m}} \ar@{->}[r]^{\psi} & \psi(G') \ar@{-}[u]_{p^{m}} &&\\
 && \psi^{-1}(G' \times \dots \times G')\ar@{-}[u] \ar@{->}[r]^{\psi}         & G' \times \dots \times
 G'  \ar@{-}[u]_{p^{m(p-1)}} &&
 }
\]
\caption{Subgroups and indices in the branching structure of $G$ ($p
\neq 2$)}\label{f:indices}
\end{figure}
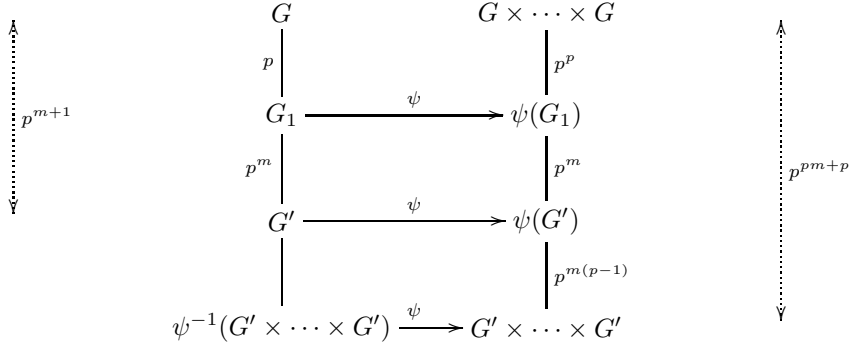
By Lemma~\ref{lemma:DB_1}, $[G:G']=p^{m+1}$. Since $[G:G_1]=p$, this
implies $[G_1:G']=p^m$, and therefore $[\psi(G_1):\psi(G')]=p^m$
($\psi$ is an embedding). Since $[G:G']=p^{m+1}$, we have $[G \times
\dots \times G:G' \times \dots \times G']=p^{pm+p}$ and therefore
\[
 [G':\psi^{-1}(G' \times \dots \times G')] = [\psi(G'):G' \times \dots \times G'] = p^{pm+p-m-p} = p^{m(p-1)}.
\]
\end{proof}

\subsection{Even prime case}

We consider now the case $p=2$ and prove Theorem~\ref{branch}(b). We
again proceed by subdividing the proof into several lemmas.

\begin{lemma}\label{l:gk}
Let $G$ be a group in $\G_{2,m}$, $m \geq 2$. Then,
\[ \overline{B_1} = B_1 \ltimes K. \]

In particular, $[\overline{B_1},K]=2^{m-1}$ and $[G:K] =2^{m+2}$.
\end{lemma}
\begin{proof}
Since $\pi_B(K)=1$ and the only element of $B_1$ with trivial
$B$-projection is $1$, we see that $B_1 \cap K = 1$. Since $K$ is
normal in $G$, it is normal in $\overline{B_1}$ as well. All that we
need to show is that every element of $\overline{B_1}$ is a product
of an element in $B_1$ and an element in $K$.

Let $g\in \overline{B_1}$. The group $G$ is generated by the
elements in $B_1$ together with $a$ and the element $d \in B_0 -B_1$
from Lemma~\ref{lemma:DB_1}. In the factor group $G/K$ all of these
generators commute except possibly $a$ and $d$. Thus the element $g$
can be written as $g=hb_1k$, for some $h \in D$, $b_1 \in B_1$ and
$k \in K$. Since $g,b_1k \in \overline{B_1}$ and $D \cap
\overline{B_1} = 1$, we must have $h=1$. Thus $g=b_1k$.

The claim on the indices then follows from $|B_1| = 2^{m-1}$.
\end{proof}

\begin{lemma}\label{overk}
Let $G$ be a group in $\G_{2,m}$, $m \geq 2$. Then $G$ is a regular
branch group over $K$.
\end{lemma}

\begin{proof}
Let $b$ be an element in $B_1$ with $\omega(b)=a$ and let $d \in
B_0$. Then $[b,a]$ is in $K$ and so is
\[ [[b,a],d] = a^{-1}b^{-1}abd^{-1}b^{-1}a^{-1}bad =
    a^{-1}b^{-1}ad^{-1}a^{-1}bad = [b^a,d]. \]
Since
\begin{multline*}
 \psi( [(b^a,d] )  = \\
    (\rho(b^{-1}),\omega(b^{-1})(\omega(d^{-1}),\rho(d^{-1}))
    (\rho(b),\omega(b))(\omega(d),\rho(d)) =
    (1,[a,\rho(d)]),
\end{multline*}
and $\rho(B_0)=B_1$ we conclude that, for $b' \in B_1$,
\[ \psi([b^a,\rho^{-1}(b')])=(1,\dots,1,[a,b']) \in \psi(K). \]
The rest of the proof proceeds just like in the case $p \neq 2$
(Lemma~\ref{overg'}).
\end{proof}

\begin{lemma}
Let $G$ be a group in $\G_{2,m}$, $m \geq 2$. Then there exists an
element $c \in B_{-1} - B_0$ and an element $d \in B_0 - B_1$ such
that
\[ c = (a,d). \]
Let
\[\hat{C} = \langle (a,d), (d,a) \rangle. \]
Then,

(a)
\[ \psi(G_1) =
 \hat{C} \ltimes (\overline{B_1} \times \overline{B_1}).
\]

(b)
\[ G \times G =
  (\hat{C} \ltimes \hat{D}) \ltimes (\overline{B_1} \times \overline{B_1}),
\]
where
\[ \hat{D} = \langle (1,a), (1,d) \rangle = 1 \times D \].

(c) the index of $\psi(G_1)$ in $G \times G$ is $2^3$.

(d) the index of $\psi^{-1}(K \times \dots \times K)$ in $K$ is
$2^m$.
\end{lemma}
\begin{proof}
The element $c$ is chosen just like in the case $p \neq 2$.

(a) The same proof as in the case $p \neq 2$.

(b) Since $G = D \ltimes \overline{B_1}$, all we need to show is
that $D \times D = \hat{C} \ltimes \hat{D}$. Note that $\hat{C}$ is
a diagonal subgroup of $D \times D$ (the first coordinate in
$\hat{C}$ is the image of the second coordinate under the
automorphism of the dihedral group $D=D_4$ defined by $a \mapsto d$,
$d \mapsto a$. Thus $D \times D = \hat{C} \ltimes (1 \times D) =
\hat{C} \ltimes \hat{D}$.

(c) Follows from (a) and (b) and the observation that
$|\hat{D}|=2^3$.

(d) All the relevant subgroups and indices are given in
Figure~\ref{f:indices2}.
\begin{figure}[!ht]
\[
\xymatrix{
 \ar@{<..>}[dd]^{2^{m+2}} &&
 G                                  & G \times G && \ar@{<..>}[ddd]^{2^{2m+4}}\\
 && G_1  \ar@{-}[u]^2 \ar@{->}[r]^{\psi}         & \psi(G_1) \ar@{-}[u]_{2^3} && \\
 && K  \ar@{-}[u]^{2^{m+1}} \ar@{->}[r]^{\psi} & \psi(K) \ar@{-}[u]_{2^{m+1}} &&\\
 && \psi^{-1}(K \times  K) \ar@{-}[u] \ar@{->}[r]^{\psi}         & K \times K
 \ar@{-}[u]_{2^{m}} &&
 }
\]
\caption{Subgroups and indices in the branching structure of $G$
($p=2$)}\label{f:indices2}
\end{figure}
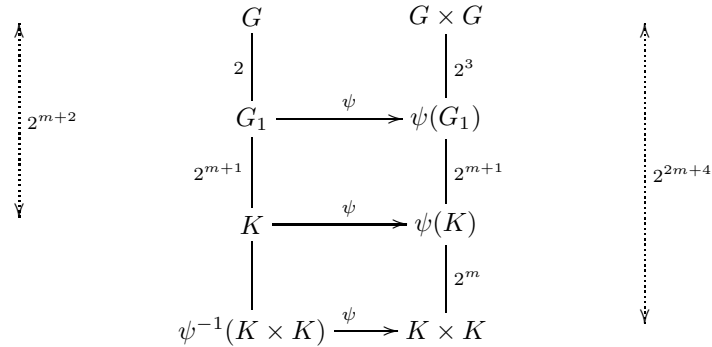

\end{proof}


\section{Hausdorff dimension}

Let $G$ be a group in $\G_{p,m}$, $m \geq 2$. Our results related to
the branching structure of $G$ allow us to calculate the Hausdorff
dimension of the closure $\overline{G}$ in the pro-$p$-group $\A(p)$
of $p$-adic automorphisms of the $p$-ary tree $X^*$. The metric used
on $\A(p)$ is the one defined in~\eqref{metric-a}. With respect to
this metric, the Hausdorff dimension of a closed subgroup
$\overline{G}$ of $\A(p)$ can be calculated
(see~\cite{barnea-s:hausdorff}) simply by comparing the relative
sizes of the level stabilizers by
\[
  \dimh(\overline{G}) = \liminf_{n \to \infty}
  \frac{\log~[\overline{G}:\overline{G}_n]}
       {\log~[\A(p):\A(p)_n]}.
\]
Since $[\A(p):\A(p)_n]= p^{\frac{p^n-1}{p-1}}$ and the group
$\overline{G}$ is the closure of $G$, the above dimension formula
reduces in our situation to
\begin{equation}\label{dimensionformula}
  \dimh(\overline{G}) = \liminf_{n \to \infty}
  \frac{p-1}{p^n-1}~\log_p~[G:G_n].
\end{equation}

\begin{theorem}\label{dimension}
Let $G$ be a group in $\G_{p,m}$, $m \geq 2$. The Hausdorff
dimension of the closure $\overline{G}$ of $G$ in $\A(p)$ is
\[
 \dimh(\overline{G}) =
  \begin{cases} 1 - \frac{1}{p^m}, & p \neq 2\\
                1 - \frac{3}{2^{m+1}}, & p = 2
 \end{cases}.
\]
In particular, the dimension $\dimh(\overline{G})$ approaches 1 as
$m$ tends to infinity.
\end{theorem}

The proof will be derived as a consequence of the following lemmas
and the dimension formula~\eqref{dimensionformula}.

\begin{lemma}\label{levelm+1}
Let $G$ be a group in $\G_{p,m}$, $m \geq 2$. Then

(a) the group $G/G_{m+1}$ is canonically isomorphic to the full
group of $p$-adic automorphisms of the finite $p$-ary tree
$X_{[m+1]}$.

(b) $[G:G_{m+1}] = p^\frac{p^{m+1}-1}{p-1}$.
\end{lemma}
\begin{proof}
(a) We identify $\pi_A(g_u)$ with the activity of $g$ at $u$ (indeed
$\pi_A(g_u)=a^k$ if and only if the activity of $g$ at $u$ is
$\pi^k$, where $\pi$ is the cyclic permutation $(0 \ \dots \ p-1)$).
A vertex $u$ with $\pi_A(g_u) \neq 1$ is called an active vertex of
$g$.

Consider the restriction of the action of the standard basis
elements $d_0,\dots,d_{m-1}$, together with $a$, on the finite tree
$X_{[m+1]}$ consisting of the top $m+1$ levels of $X^*$, i.e.,
consider their activity at all vertices up to level $m$. Denote the
vertex $(p-1)^{i-1}0$ at level $i$, $i=1,\dots,m$, by $v_i$ and
denote the root by $v_0$. The only active vertex of $d_0$ is $v_m$,
where the activity is $a$. For $d_1$, the vertex $v_{m-1}$ is
active, with activity $a$, and the only other possibly active vertex
is $v_m$. In general, for $d_i$, $i=0,\dots,m-1$, the activity at
$v_{m-i}$ is $a$ and the only other possibly active vertices are
$v_{i+1},\dots,v_m$. By changing the generators $d_i$,
$i=0,\dots,m-1$, to $d_i'$, $i=0,\dots,m-1$, where
\begin{align*}
 d_0' &= d_0 \\
 d_1' &= d_0^*d_1 \\
 d_2' &= d_0^*d_1^*d_2 \\
 \dots \\
 d_{m-1}' &= d_0^*d_1^* \dots d_{m-2}^* d_{m-1},
\end{align*}
and the $*$'s are appropriately chosen integers, we obtain a
generating set $d_i'$, $i=0,\dots,m-1$, for $B$ such that the only
active vertex of $d_i'$ is $v_i$, where the activity is $a$.
Together with $a$, whose only activity occurs at $v_0$ these tree
automorphisms act on $X_{[m+1]}$ as the full group od $p$-adic
automorphisms.

(b) From (a) $|G/G_{m+1}| = |\underbrace{C_p \wr C_p \wr \dots \wr
C_p}_{m+1}| = \frac{p^{m+1}-1}{p-1}$.
\end{proof}

\begin{lemma}
Let $G$ be a group in $\G_{p,m}$, $m \geq 2$.

(a) the stabilizer $G_{m+1}$ is contained in the commutator $G'$.

(b) when $p=2$, the stabilizer $G_{m+1}$ is contained in $K$.
\end{lemma}
\begin{proof}
(a) Let $g \in G$. If $\pi_A(g) \neq 1$, then $g$ is active at the
root. Consider $\pi_B(g)$ and write it in terms of the generating
set $d_i'$, $i=0,\dots,m-1$. If $\pi_B(g) \neq 1$, then one of these
generators, say $d_i'$, will appear with an exponent $k$ from
$\{1,\dots,p-1\}$. But this means that there must be some active
vertices of $g$ at level $m-i$ (in fact, $a^k$ is the product of all
activites of $g$ at level $m-i$). Thus if $g \not \in G'$, then
there is some activity on at least one of the levels $0,\dots,m$.
Our claim follows by contraposition.

(b) Let $d_i$, $i=0,\dots,m-1$, be the standard basis elements of
$B$ and let $f(x)=x^m + a_{m-1}x^{m-1} + \dots + a_1x + a_0$ be the
polynomial defining $G$. Then,
\begin{align*}
 [a,d_i ] &= (d_{i+1},d_{i+1}), \text{ for } i =0,\dots,m-2 \\
 [a,d_{m-1}] &= (ca,ac),
\end{align*}
where $c = d_0^{a_0}d_1^{a_1} \dots d_{m-1}^{a_{m-1}}$. In
particular, the order of $[a,d_0]$ is 2.

Denote the commutator $[a,d_i]$, $i=0,\dots,m-1$, by $\delta_i$.
From part (a) we already know that $G_{m+1} \leq G'$. From
Lemma~\ref{lemma:DB_1} and Lemma~\ref{l:gk} it follows that $G = D
\ltimes (B_1 \ltimes K)$, where $D=\langle a, d_0 \rangle$. In
particular, $G=DB_1K$. Since $K \leq G'$ and the only elements in
$DB_1$ that are in $G'$ are $1$ and $\delta_0$ we have $G' = \langle
\delta_0 \rangle \ltimes K = C_2 \ltimes K$.

Let $g \in G_{m+1}$ and let it be written as a product of conjugates
of $\delta_i$, $i=0,\dots,m-1$ (this is possible since the first
$m-1$ have order 2 and the last one is inverted by conjugating it by
$a$). Consider now the section $g_0=\varphi_0(g)$ of $g$. The
section $g_0$ is in $G_ m$ and is an element of $\langle
d_1,\dots,d_{m-1},ca \rangle^G$. Since the projections $\pi_A$ and
$\pi_B$ are constant on conjugacy classes of $G$, their images
depend only on the parity of the number of conjugates of $d_i$,
$i=1,\dots,m-1$, and $ca$ used to express $g_0$, which in turn
depend only on the parity of the number of conjugates of $\delta_i$,
$i=0,\dots,m-1$ used to express $g$. Since $\pi_A(g_0)$ must be 1,
we see that the number of conjugates of $ca$ used in $g_0$ is even.
Thus $\pi_B(g_0)$ depends only on the number of conjugates of $d_i$,
$i=1,\dots,m-1$, and is an element in $B_1=\langle d_1, \dots,
d_{m-1}\rangle$. We consider the action of $g_0$ on level $m$ (it
has to be trivial) and the activity only up to level $m-1$. There is
a change of generators of $B_1$,
\begin{align}
 d_1' &= d_1 \notag \\
 d_2' &= d_1^*d_2 \label{d'd2}\\
 \dots \notag\\
 d_{m-1} &= d_1^*d_2^* \dots d_{m-1}^*d_{m-1},  \notag
\end{align}
where the $*$'s represents appropriately chosen integers, such that
$d_i'$, $i=1,\dots,m-1$, has only one active vertex (up to level
$m-1$), which is at level $m-i$. Rewrite $g_0$ in terms of
conjugates of $d_i'$. Since conjugation does no change the parity of
the total activity at some level, we see that in order for $g_0$ to
stabilize level $m$, the parity of the number of conjugates of
$d_i'$ used to express $g_0$ must be even, for $i=1,\dots,m-1$.
Using \eqref{d'd2} to express $g_0$ back in terms of conjugates of
$d_i$, we see that the number of conjugates of $d_i$,
$i=1,\dots,m-1$, as well as the number of conjugates of $ca$, used
in $g_0$ must be even, and therefore the number of conjugates of
$\delta_i$, $i=0,\dots,m-1$  used in $g$ must be even. In
particular, $g$ must be in $K$.
\end{proof}

The following simple lemma applies to all regular branch groups (not
necessarily subgroups of $\A(p)$).

\begin{lemma}\label{stabstab}
Let $G$ be a regular branch group acting on the $k$-ary tree $X^*$,
branching over a subgroup $K$ that contains the level stabilizer
$G_s$. Then, for $n \geq s$,

\[ \psi(G_{n+1}) = G_n \times \dots \times G_n. \]
\end{lemma}
\begin{proof}
The inclusion $\psi(G_{n+1}) \leq G_n \times \dots \times G_n$ is
valid for any group of tree automorphisms. For the converse, if
$(g_0,\dots,g_{k-1}) \in G_n \times \dots \times G_n$, then
$(g_0,\dots,g_{p-1}) \in K \times \dots \times K$ and therefore
there exists an element $g \in G$ (in fact in $K$) such that
$\psi(g) = (g_0,\dots,g_{k-1})$. The element $g$ must come from
$G_{n+1}$ (otherwise at least one of its sections would not be in
$G_n$).
\end{proof}

\begin{proposition}\label{p:dimension}
Let $G$ be an infinite self-similar group of $p$-adic automorphisms.
Further, for some $s \geq 1$ let $[G:G_s]=p^{\frac{r}{p-1}}$, let
$[G \times \dots \times G:\psi(G_1)]=p^t$, and let $\psi(G_{n+1}) =
G_n \times \dots \times G_n$, for $n \geq s$. Then

(a) for $n \geq s$,
\[ [G:G_n] = p^{\frac{r-t+1}{p-1}p^{n-s}+\frac{t-1}{p-1}}.\]

(b) the Hausdorff dimension of $\overline{G}$ is
\[ \dimh(\overline{G}) = \frac{r-t+1}{p^s} . \]
\end{proposition}
\begin{proof}
(a) We first note that $[G:G_1]=p$. Indeed, the group $G$ is not
trivial and therefore there is a nontrivial activity at some vertex.
By self-similarity, there must be nontrivial activity at the root as
well.

We proceed by induction on $n$. The formula is clearly correct for
$n =s$. For $n > s$, we have
\begin{align*}
 [G:G_n] &= [G:G_1][G_1:G_n] = p~[\psi(G_1):\psi(G_n)] = \\
  &= p~\frac{[G \times \dots \times G:\psi(G_n)]}
            {[G \times \dots \times G:\psi(G_{1})]} =
     p~\frac{[G \times \dots \times G:G_{n-1} \times \dots \times G_{n-1}]}
         {[G \times \dots \times G:\psi(G_{1})]} = \\
    &= p^{1-t}~[G:G_{n-1}]^p = p^{\frac{r-t+1}{p-1}p^{n-s}+\frac{t-1}{p-1}}.
\end{align*}

(b)
\begin{align*}
 \dimh(\overline{G}) &=
 \liminf_{n \to \infty} \frac{p-1}{p^n-1}~\log_p~[G:G_n] = \\
 &= \liminf_{n \to \infty} \frac{p-1}{p^n-1}~\left(\frac{r-t+1}{p-1}p^{n-s}+\frac{t-1}{p-1}\right) =
 \frac{r-t+1}{p^s}.
\end{align*}
\end{proof}

\begin{corollary}
Let $G$ be a group in $\G_{p,m}$, $m \geq 2$. Then, for $n \geq
m+1$,
\[ [G:G_n] = p^{\frac{p^n - tp^{n-m-1}+t-1}{p-1}},\]
where $t=p$, when $p \neq 2$, and $t=3$, when $p=2$.
\end{corollary}
\begin{proof}
Direct corollary of Proposition~\ref{p:dimension}. Note that in our
case $s=m+1$, $r=p^{m+1}-1$ and $t=p$, when $p \neq 2$, and $t=3$,
when $p=2$.
\end{proof}

\begin{proof}[Proof of Theorem~\ref{dimension}]
We again use Proposition~\ref{p:dimension}. We have $s=m+1$,
$r=p^{m+1}-1$ and $t=p$, when $p\neq 2$, and $t=3$, when $p=2$. Thus
\[ \dimh(\overline{G}) = \frac{r-t+1}{p^s} =
\frac{p^{m+1}-t}{p^{m+1}} = 1- \frac{t}{p^{m+1}}.
\]
\end{proof}


\section{Finitely constrained groups}

Let $\mathcal P$ be a subgroup of the group $\Aut(X_{[s]})$ of
automorphisms of the finite $k$-ary tree of size $s$. We call
$\mathcal P$ a group of patterns (of size $s$) and its complement
$\mathcal F$ in $\Aut(X_{[s]})$ a set of forbidden patterns. An
automorphism $f$ of $X^*$ does not contain a forbidden pattern from
$\F$ if there is no section of $f$ whose action on $X^s$ agrees with
the action of some element in $\F$ (in other words, there is no
section of $f$ whose activity at each vertex on levels 0 through
$s-1$ agrees with the activity of some element in $\F$).

\begin{definition}
A group $G$ of tree automorphisms of $X^*$ is a finitely constrained
group if there exists a finite set of forbidden patterns $\F$ of
size $s$ such that $G$ consists of all elements in $\Aut(X^*)$ that
do not contain forbidden patterns from $\F$.
\end{definition}

The following proposition is proved in~\cite{grigorchuk:unsolved}

\begin{proposition}[\cite{grigorchuk:unsolved}, Proposition 7.5]\label{slava}
Let $G$ be a spherically transitive, finitely constrained group of
tree automorphisms of $X^*$. Then $G$ is a self-similar, regular
branch group.
\end{proposition}

We will provide a characterization of spherically transitive,
finitely constrained groups in terms of branch groups. Recall that a
congruence subgroup of a group of tree automorphisms is any subgroup
that contains a level stabilizer.

\begin{theorem}\label{characterization}
Let $G$ be a group of tree automorphisms of $X^*$. The following are
equivalent:

(i) $G$ is a spherically transitive, finitely constrained group
(with forbidden patterns of size $s+1$).

(ii) $G$ is the closure of a self-similar, regular branch group $H$,
branching over a congruence subgroup (containing the stabilizer
$H_s$ of level $s$).

(iii) $G$ is the closure of a self-similar regular branch group $H$,
branching over a level stabilizer (level $s$).
\end{theorem}
\begin{proof}
(i) implies (iii). Evident from the proof of Proposition~7.5
in~\cite{grigorchuk:unsolved} (the proof itself claims that $G$ is
the closure of a self-similar regular branch group $H$, branching
over $H_{s+1}$, but a careful reading reveals that the argument
works equally well for $H_s$).

(iii) implies (ii). By definition.

(ii) implies (i). Let $G$ be the closure of a self-similar, regular
branch group $H$, branching over a congruence subgroup $K$ that
contains $H_s$.

Since $H$ is spherically transitive, so is $G$. Also, note that $G$
itself is self-similar (each section $g_u$, in the limit $g=\lim
h_n$ of a sequence of automorphisms $\{h_n\}$ in $H$, is the limit
$g_u=\lim (h_n)_u$ of the sequence of corresponding sections, which
are all in $H$, by the self-similarity of H).

It follows from Lemma~\ref{stabstab} that, for $n \geq s$,
$\psi(H_{n+1}) = H_n \times \dots \times H_n$.

Let $\F$ be the complement in $\Aut(X_{[s+1]})$ of $H/H_{s+1}$ and
let $G(\F)$ be the finitely constrained group of $p$-adic
automorphisms defined by the forbidden patterns in $\F$. If $g$ in
$G$ then $g$ cannot contain a pattern from $\F$. Indeed, if it did,
then one of its sections, which by self-similarity is also in $G$,
would agree with an element in $\F$ on $X^{s+1}$. However, this is
impossible, since this would mean that an element in $H$ agrees with
a forbidden pattern in $\F$ on $X^{s+1}$. Thus $G \leq G(\F)$.

Let $g$ be an element in $G(\F)$. We will prove by induction on $n$
that $g = \lim h_n$, for some sequence $\{h_n\}$ of elements in $H$.
In fact, we will show that, for every $n \geq 0$, there exists an
element $h_n$ in $H$ such that $h_n^{-1}g \in \Aut(X^*)_n$ (recall
that $\Aut(X^*)_n$ is the stabilizer of level $n$ in $ \Aut(X^*)$).
This is certainly true for $n=s+1$. Indeed, $g \in G(\F)$ and
therefore its action on $X^n$ agrees with the action of some element
$h_{s+1} \in H$. Assume $n \geq s+2$ and that the inductive claim is
true for values smaller than $n$. Multiply $g$ by $h_{n-1}$ to
obtain an element $f=h_{n-1}^{-1}g \in \Aut(X^*)_{n-1}$. This
element is still in $G(\F)$ and so are its sections,
$f_0,\dots,f_{k-1}$. By inductive hypothesis, there are elements
$f_0',\dots,f_{k-1}'$ in $H$ such that $(f_i')^{-1}f_i \in
\Aut(X^*)_{n-1}$, $i=0,\dots,k-1$. The elements $f_i'$,
$i=0,\dots,k-1$, are elements in $H_{n-2}$, and therefore there
exists an element $f'$ in $H_{n-1}$ such that
$f'=(f_0',\dots,f_{k-1}')$. Then $(f')^{-1}f =
((f_0')^{-1}f_0,\dots,(f_{k-1}')^{-1}f_{k-1}) \in \Aut(X^*)_n$. Thus
$\Aut(X^*)_n \ni (f')^{-1}f = (f')^{-1}h_{n-1}^{-1}g$ and we may
take $h_n = h_{n-1}f'$.
\end{proof}

\begin{corollary}
Let $G$ be a group in $\G_{p,m}$. $m \geq 2$. Then the closure
$\overline{G}$ is a finitely constrained subgroup of $\A(p)$ with
forbidden patterns of size $m+2$.
\end{corollary}

\begin{theorem}\label{generaldimension}
(a) Any infinite, finitely constrained group of $p$-adic
automorphisms has positive and rational Hausdorff dimension.

(b) Let $G$ be an infinite finitely constrained subgroup of the
iterated wreath product $\A_Q = Q \wr Q \wr \dots$ of subgroups of
$\Sym_k$ and let $G$ be the closure of the self-similar, regular
branch group $H$, branching over the level stabilizer $H_s$. Then
the Hausdorff dimension of $G$ (with respect to $\A_Q$) is positive
and it is equal to
\[ \dimh(G) = \frac{r-t+\epsilon}{k^s}, \]
where $k$ is the arity of the tree, $q$ is the order of $Q$,
$[H:H_1]=q^\epsilon$, $[H\times \dots \times H:\psi(H_1)]=q^t$ and
$[H:H_s]=q^{\frac{r}{k-1}}$.
\end{theorem}
\begin{proof}
(a) By Theorem~\ref{characterization} and Lemma~\ref{stabstab},
Proposition~\ref{p:dimension} applies. Note that any infinite
self-similar group of $p$-adic automorphisms is transitive (by an
unpublished argument of D.~Savchuk).

Since $G$ is infinite, the sequence of indices $[G:G_n]$ is
unbounded. It follows from Proposition~\ref{p:dimension}(a) that
$r-t+1$ must be strictly positive, which, by
Proposition~\ref{p:dimension}(b), makes the Hausdorff dimension of
$\overline{G}$ strictly positive. The dimension is rational since
$p$, $r$, $s$ and $t$ are integers.

(b) Following the exact same line of reasoning as in the proof of
Proposition~\ref{p:dimension} one can calculate the Hausdorff
dimension explicitly and conclude that it must be positive. Indeed,
for $n \geq s$, one gets the recursive formula $[H:H_{n+1}] =
q^{\epsilon-t}[H:H_n]^k$, with initial condition
$[H:H_s]=q^{\frac{r}{k-1}}$, whose unique solution, for $n \geq s$,
is
\[
 [H:H_n] = q^{\frac{r-t+\epsilon}{k-1}k^{n-s} + \frac{t-\epsilon}{k-1}}
\]
and the claims easily follow.
\end{proof}

Theorem~\ref{generaldimension}(a) answers positively
Problem~7.1.(ii) on rationality of Hausdorff dimension in the
context of infinite, finitely constrained groups of $p$-adic
automorphisms.

However, the rationality of $\dimh(G)$ is not guaranteed in
Theorem~\ref{generaldimension}(b), since $r$, $t$, and $\epsilon$
are not necessarily integers. The following example shows that the
dimension can be irrational even for closures of bounded automaton
groups.

\begin{example}
It is known that the closure $G$ of the Hanoi group $H$ on 3 pegs
(see~\cite{grigorchuk-s:hanoi-cr} for a definition of $H$ as an
automaton group) is finitely constrained
group~\cite{grigorchuk-n-s:oberwolfach2} with forbidden patterns of
size $2$. Thus $G$ is a regular branch group over its level
stabilizer $G_1$. Since $[G:G_1]=[H:H_1]=6$ and $[G \times G \times
G:\psi(G_1)]=2$ we get $\epsilon = 1$, $r=2$ and $t=\log_6(2)$, and
therefore, by Theorem~\ref{generaldimension}(b) the Hausdorff
dimension is
\[ \dimh(G) = 1 - \frac{1}{3}\log_6(2) \approx 0.871. \]

Note that it would be wrong to use $[H\times H \times
H:\psi(H_1)]=2^5$ instead of $[G \times G \times G:\psi(G_1)]=2$ in
the above calculation. The difference comes from the fact that, even
though $G$ is a finitely constrained group that is the closure of
the self-similar, regular branch group $H$, $H$ is not branching
over a congruence subgroup and therefore we cannot use the index
$[H\times H \times H:\psi(H_1)]=2^5$.
\end{example}

Regarding Problem~7.1.(ii) from~\cite{grigorchuk:unsolved}, it would
still be interesting to check if the Hausdorff dimension of closed,
self-similar groups of $p$-adic automorphisms is always rational.

\section{Further finiteness properties}

We quickly list some further properties of the groups in $\G$, some
of which come as consequences of known facts about some larger
classes of groups that include $\G$.

\begin{proposition}
Every group in $\G$ is amenable.
\end{proposition}
\begin{proof}
This follows from the fact that each group in $\G$ is defined by a
bounded automaton and the general result on amenability of such
groups in~\cite{bartholdi-k-n-v:bounded}.
\end{proof}

\begin{proposition}\label{torsion}
Let $G$ be a group in $\G_{p,m}$, $m \geq 2$. The following are
equivalent:

(i) $G$ is a $p$-group.

(ii) there exists $r$ such that
\[ B_0 \cup B_1 \cup \dots \cup B_{r-1} = B.  \]

(iii) every nontrivial $\rho$-orbit intersects the kernel
$B_0=\Ker(\omega)$.
\end{proposition}
\begin{proof}
(ii) is equivalent to (iii). Clear.

(ii) implies (i). Follows from general results on spinal groups
in~\cite{bartholdi-s:wpg}.

(i) implies (iii). Assume there is a $\rho$-orbit $b_0 \mapsto b_1
\mapsto \dots \mapsto b_{s-1} \mapsto b_0$, $s \geq 1$, that is
completely outside of $B_0$. Then, for $i=0,\dots,s-1$,
\[ b_i = (a^{n_{i+1}},\dots,b_{i+1}), \]
with indices taken modulo $s$, where each of the numbers $n_i$ is in
$\{1,\dots,p-1\}$.

Assume that $G$ is a $p$-group. Denote the order of an element $g$
by $\Omega(g)$. Consider $a^{n_i}b_i$. It does not fix the first
level and therefore it is not trivial. Since each of the first level
sections of $(a^{n_i}b_i)^p$ is conjugate to $a^{n_{i+1}}b_{i+1}$ we
have that $\Omega(a^{n_i}b_i) = p\cdot\Omega(a^{n_{i+1}}b_{i+1})$.
Therefore, $\Omega(a^{n_0}b_0) = p^s\cdot\Omega(a^{n_0}b_0)$, which
is a contradiction.
\end{proof}

In fact, it seems that the conditions in Proposition~\ref{torsion}
are equivalent to the condition that the polynomial $f$ defining
$G_{p,f}$ does not have a binomial factor (of any degree). This is
clearly always true in one direction (namely binomial factors induce
elements of infinite order). The converse is also relatively
straightforward for $p=2$.

\begin{proposition}
Let $G$ be a $p$-group in $\G_{p,m}$, $m \geq 2$, and let $r$ be the
smallest integer such that $B_0 \cup B_1 \cup \dots \cup B_{r-1} =
B$. Then

(a) then there exists a positive constant $C$ such that the order of
every element of $G$ of length $n$ is at most
\[ Cn^{(r-1)\log_2 p}, \]

(b) $G$ has intermediate growth and there exists positive constants
$c_1$ and $c_2$ such that the growth function $\gamma$ of $G$
satisfies the inequalities
\[ c_1^{\sqrt{n}} \leq \gamma(n) \leq c_2^{n^\alpha}, \]
where $\alpha=\frac{\log p}{\log p - \log \eta}$ and $\eta$ is the
positive root of the polynomial $x^r + x^{r-1} + x^{r-2} - 2$.
\end{proposition}
\begin{proof}
Follows from general results on spinal groups
in~\cite{bartholdi-s:wpg}.
\end{proof}

Note that if every nontrivial $\rho$-orbit of $B$ intersects $B_0$
nontrivially, then $r$ is the smallest integer such that $B_0 \cup
\dots \cup B_{r-1} = B$ if and only if $r-1$ is the largest length
of a piece of a $\rho$-orbit that lies entirely outside of $B_0$. If
$f$ is a primitive polynomial of degree $m$ over the prime $p$, then
there is only one nontrivial $\rho$-orbit and the corresponding
value of $r$ is $m+1$. In a sense, the examples that are closest to
the first Grigorchuk group are the ones defined by primitive
polynomials over $p$.

\def\cprime{$'$}



\end{document}